\documentclass[twoside,leqno]{article}

\usepackage[letterpaper]{geometry}

\usepackage{siamproceedings}

\usepackage[T1]{fontenc}
\usepackage{amsfonts}
\usepackage{euscript}
\usepackage{amssymb}
\usepackage{graphicx}
\usepackage{epstopdf}
\usepackage[shortlabels]{enumitem}
\usepackage{algorithmic}
\ifpdf
  \DeclareGraphicsExtensions{.eps,.pdf,.png,.jpg}
\else
  \DeclareGraphicsExtensions{.eps}
\fi

\newsiamremark{remark}{Remark}
\newsiamremark{hypothesis}{Hypothesis}
\crefname{hypothesis}{Hypothesis}{Hypotheses}
\newsiamthm{claim}{Claim}
\newsiamremark{example}{Example}

\usepackage{hyperref}
\usepackage{latexsym}
\usepackage{todonotes}
\usepackage{nicefrac}
\usepackage{epsfig}
\usepackage{stmaryrd}
\usepackage{setspace}
\usepackage{tikz-cd}

\usetikzlibrary{decorations.pathreplacing} 

\ifpdf
\usepackage{pdfsync}
\fi
\usetikzlibrary{positioning,matrix,arrows,decorations.pathmorphing}
\usepackage{enumitem}

\usetikzlibrary{positioning,arrows}
\usetikzlibrary{decorations.markings}
\usepackage{caption, subcaption}
\usetikzlibrary{arrows.meta,
                bbox}

\renewcommand{\Im}{\operatorname{Im}}

\renewcommand{\dim}{\operatorname{dim}}

\newcommand{\cP}{\mathcal{P}}

\newcommand{\cR}{\mathcal{R}}

\newcommand{\C}{{\mathbb{C}}}
\newcommand{\Z}{{\mathbb{Z}}}
\newcommand{\N}{{\mathbb{N}}}

\newcommand{\R}{{\mathbb{R}}}

\newcommand{\IC}{\operatorname{IC}}
\newcommand{\IH}{\operatorname{IH}}

\renewcommand{\a}{\alpha}

\newcommand{\cs}{\mathbb{C}^\times}

\newcommand{\cA}{\mathcal{A}}

\newcommand{\cE}{\mathcal{E}}

\newcommand{\cF}{\mathcal{F}}
\newcommand{\cO}{\mathcal{O}}

\newcommand{\cG}{\mathcal{G}}
\newcommand{\cH}{\mathcal{H}}
\newcommand{\cM}{\mathcal{M}}

\newcommand{\cL}{\mathcal{L}}

\renewcommand{\and}{\qquad\text{and}\qquad}

\newcommand{\Sym}{\operatorname{Sym}}

\renewcommand{\H}{{\operatorname{H}}}
\newcommand{\CH}{{\operatorname{CH}}}

\renewcommand{\P}{\mathbb{P}}

\newcommand{\beq}{\begin{eqnarray*}}
\newcommand{\eeq}{\end{eqnarray*}}
\newcommand{\rk}{\operatorname{rk}}

\usepackage{todonotes}

\newcommand{\excise}[1]{}
\usepackage{amsopn}

\begin{document}

\title{\Large Intersection cohomology without spaces}
    \author{Tom Braden\thanks{Department of Mathematics and Statistics, University of Massachusetts, Amherst, MA (\email{braden@umass.edu}).}
    \and Nicholas Proudfoot\thanks{Department of Mathematics, University of Oregon, Eugene, OR
  (\email{njp@uoregon.edu}).}}

\date{}

\maketitle


\fancyfoot[R]{\scriptsize{Copyright \textcopyright\ 20XX by SIAM\\
Unauthorized reproduction of this article is prohibited}}





\begin{abstract} 
We survey three settings in which dimensions of intersection cohomology groups of algebraic varieties provide deep combinatorial and representation-theoretic information, and  computations of the groups themselves have been made using combinatorial sheaves on finite posets.  These settings are (1) intersection cohomology of Schubert varieties, the associated Kazhdan--Lusztig polynomials and their realizations via moment graph sheaves and Soergel bimodules; (2) intersection cohomology of toric varieties, the associated $g$-polynomials of convex polytopes, and their realization via the theory of intersection cohomology of fans; and (3) intersection cohomology of arrangement Schubert varieties, the associated Kazhdan--Lusztig polynomials of matroids, and their realization via intersection cohomology of matroids.  In all three settings these constructions are valid in more general situations where the variety does not exist, leading to ``intersection cohomology without spaces.''	
We give parallel presentations of these three stories, highlighting applications to KLS-polynomials.
\end{abstract}

\section{Introduction.}\label{sec:intro}
Given a finite ranked poset $P$ along with an additional piece of data called a $P$-kernel, one may assign
to each interval $x\leq y$ in $P$ a {\bf KLS-polynomial} $f_{xy}(t)\in\Z[t]$.
Examples include Kazhdan--Lusztig polynomials of Coxeter groups, $g$-polynomials of polytopes, and Kazhdan--Lusztig polynomials of matroids.
In each of these three cases, the polynomials are known to have non-negative coefficients, and in each case the proof followed
a similar historical arc.  

Non-negativity was first proved in special cases by interpreting the polynomials as local intersection cohomology Poincar\'e polynomials
for certain algebraic varieties.  The special class of Coxeter groups is the class of Weyl groups, and the related geometry
is that of flag varieties and their Schubert stratifications \cite{KL80}.  
The special polytopes are the rational ones, which correspond to toric varieties \cite{DL,Fieseler}.
Finally, the special matroids are those that are realizable by hyperplane arrangements, 
and the associated geometric objects are called arrangement Schubert varieties \cite{EPW}.
A unified treatment of the geometric proofs of non-negativity in these three cases appears in \cite{KLS}.

For general Coxeter groups, polytopes, and matroids, there are no analogous spaces whose intersection cohomology groups provide geometric interpretations
of these polynomials.  However, in each case, one can give a purely algebraic construction of graded vector spaces whose Poincar\'e polynomials are
the KLS-polynomials.  These vector spaces are (subquotients of) Soergel bimodules \cite{Soe90}, intersection cohomology of fans \cite{BBFK,BBFK2}, and (quotients of) intersection cohomology modules of matroids \cite{SHCG}.  In the special cases in which there is a geometric object, these vector spaces are canonically isomorphic to the corresponding intersection cohomology groups.  In general, however, we may regard them as ``intersection cohomology without spaces''.

Our aim is to review the constructions of these three purely algebraic objects in as unified a manner as possible, stressing both the parallels and
the differences.  This review will begin with a procedure for {\em computing} the relevant intersection cohomology groups when there is a variety,
which can then be translated into a procedure for {\em defining} these groups when there is no variety.
We note that while constructing the objects involves elementary combinatorial algebra,  proving that their Poincar\'e polynomials coincide
with the corresponding KLS-polynomials is much harder. The proofs in all three cases require establishing analogues
of the hard Lefschetz theorem and the Hodge--Riemann bilinear relations for intersection cohomology of projective varieties.
In the case of Coxeter groups, 24 years elapsed between the introduction of Soergel bimodules and the resolution
of the Kazhdan--Lusztig conjectures \cite{EW14}.  The analogous gap for polytopes was only five years \cite{Karu}, and for matroids it was nonexistent.
However, even in the matroid case, constructing the objects is much easier than computing their dimensions.  While we will say a few words about the proofs
of the Hodge-theoretic statements, our primary focus will be on the constructions of the vector spaces.

\section{Intersection cohomology.}
The cohomology ring of a smooth projective algebraic variety $X$ over the complex numbers satisfies Poincar\'e duality, the hard Lefschetz
theorem, and the Hodge--Riemann bilinear relations.  These results also hold when $X$ is projective but singular, provided that cohomology
is replaced by (middle perversity) intersection cohomology, which is a graded module over cohomology.

The intersection cohomology of $X$ is by definition the cohomology of the intersection complex $\IC_X$, 
which is a constructible complex of sheaves on $X$, just as the de Rham cohomology of a manifold is cohomology of the complex of differential forms.  
Unlike the de Rham complex, which is locally uniform, the intersection complex has different local properties at different points, 
depending on how singular they are.  So, in addition to the global intersection cohomology $\IH(X) := \H(\IC_X)$, 
one is often also interested in the local groups $\IH(X)_p := \H(\IC_X\!|_p)$ at points $p\in X$, 
or more generally $\IH(X)_A := \H(\IC_{X}\!|_A)$ for a locally closed subspace $A\subset X$.\footnote{All of these intersection cohomology
groups are graded, but we typically suppress the grading index to minimize notation.  Thus $\IH(X)$ means the same thing as $\IH^*(X)$.}

When $X$ is smooth, $\IC_X$ is isomorphic to a constant sheaf, $\IH(X)$ is isomorphic to $\H(X)$, and all local intersection cohomology
groups are 1-dimensional and concentrated in degree zero.
For general singular varieties computing intersection cohomology requires complicated sheaf and spectral sequence calculations, but  we will see that for some nice classes of varieties elementary algorithmic computations are possible.

There are several things one could mean when one asks to compute intersection cohomology.  For instance, one could ask for:
\begin{itemize}
	\item formulas for the dimensions of $\IH(X)$ and dimensions of $\IH(X)_p$ for various points $p$, 	
	\item a construction of graded vector spaces that are canonically isomorphic to $\IH(X)$ or $\IH(X)_p$, or
	\item a construction of $\IH(X)$ as a graded module over the ordinary cohomology ring $\H(X)$.   
\end{itemize}
Dimension formulas for intersection cohomology groups, usually organized as Poincar\'e polynomials, have important applications to the combinatorics of polytopes, hyperplane arrangements, and matroids, as well as to the representation theory of Coxeter groups, Hecke algebras, Lie algebras, and algebraic groups.  One fundamental property of these polynomials is that they have non-negative coefficients, which in general cannot be seen without knowing that they are dimensions of vector spaces.  

Having canonical descriptions of the underlying vector spaces can also be very useful.  For instance, there are often natural maps between intersection cohomology groups of different spaces, and showing that these maps are injective or surjective can lead to interesting inequalities among the polynomials.  In addition, when $X$ is acted on by a group of symmetries, that action passes to intersection cohomology, thus enriching the polynomials with integer coefficients to polynomials with coefficients in the representation ring of the symmetry group.

Knowledge of $\IH(X)$ as a module over $\H(X)$ provides additional structure which often facilitates explicit calculation.  (As we discuss later, this is even more true when one considers equivariant intersection cohomology for an action of a torus.)  For instance, it is often possible to recover all of the local groups $\IH(X)_p$ from the module $\IH(X)$.  Furthermore, the cohomology functor from 
the derived category $D^b(X)$ to $\H(X)$-modules is frequently full and faithful on intersection complexes of subvarieties $Y \subset X$.  That implies that the only grading preserving endomorphisms of $\IH(Y)$ as an $\H(X)$-module are multiplication by scalars, and it allows one 
to use homomorphisms between such modules to study the abelian category 
of perverse sheaves $\cP(X) \subset D^b(X)$. 

\section{KLS-polynomials.}\label{sec:KLS}
Suppose that $P$ is a poset equipped with a strictly increasing rank function $\rk:P\to\Z$.
For any $x\leq y$, we will write $r_{xy} := \rk(y)-\rk(x)$.
A collection of polynomials $\{\kappa_{xy}(t)\mid x\leq y\}$ is called a {\bf \boldmath{$P$}-kernel} if the following conditions hold:
\begin{itemize}
\item For all $x\in P$, $\kappa_{xx}(t) = 1$.
\item For all $x\leq y\in P$, $\deg\kappa_{xy}(t)\leq r_{xy}$.
\item For all $x< z\in P$, $\displaystyle \sum_{x\leq y\leq z}t^{r_{xy}}\kappa_{xy}(t^{-1})\kappa_{yz}(t) = 0$.
\end{itemize}
Given a $P$-kernel, there is a unique collection of polynomials $\{f_{xy}(t)\mid x\leq y\}$
satisfying the following conditions:
\begin{itemize}
\item For all $x\in P$, $f_{xx}(t) = 1$.
\item For all $x < y\in P$, $\deg f_{xy}(t)< r_{xy}/2$.
\item For all $x\leq z\in P$,\;\; $\displaystyle t^{r_{xz}}f_{xz}(t^{-1}) = \sum_{x\leq y\leq z}\kappa_{xy}(t)f_{yz}(t)$.
\end{itemize}
The polynomials $f_{xy}(t)$ are called {\bf KLS-polynomials}, named for Kazhdan and Lusztig, 
who considered the special case of Example \ref{ex:coxeter} \cite{KL79}, and for Stanley, who treated
the general case \cite{Stanley-h}. 
It is not obvious that the $P$-kernel condition is precisely what is needed to guarantee the existence of KLS-polynomials, but it is true!

\begin{example}\label{ex:coxeter}
Let $W$ be a Coxeter group, equipped with its Bruhat order and ranked by the length function $\ell:W\to\N$.  
Kazhdan and Lusztig recursively defined a collection of polynomials $\{R_{vw}(t)\mid v\leq w\in W\}$ that form a $W$-kernel,
and the associated polynomials $f_{vw}(t)$ are called {\bf Kazhdan--Lusztig polynomials}.
These polynomials appear as entries of transition matrices relating two natural bases for the Hecke algebra, and (in the case of Weyl groups)
in formulas relating the characters of simple modules and Verma modules for the corresponding Lie algebra.

Suppose that $W$ is the Weyl group associated with a reductive algebraic group $G$, and consider the Schubert stratification 
of the associated flag variety by orbits of the Borel subgroup:
$G/B = \bigsqcup X_w$.  
The dimension of $X_w$ is equal to $\ell(w)$, and $X_v\subset\overline{X}_w$ if and only if $v\leq w$.  Kazhdan and Lusztig \cite{KL80} proved
that $f_{vw}(t)$ coincides with the Poincar\'e polynomial for the local intersection cohomology\footnote{All cohomology groups
will be implicitly taken with coefficients in $\R$ until Section \ref{sec:positive}.} of $\overline{X_w}$
at a point $p\in X_v$:
\begin{equation}\label{Weyl-IH}f_{vw}(t) = \sum_{i\geq 0} t^i \dim \IH^{2i}(\overline{X_w})_p.
\end{equation}
In particular, this implies that the coefficients are non-negative.  Non-negativity for arbitrary Coxeter groups was an open conjecture until it was proved by 
Elias and Williamson \cite{EW14}, using the theory of Soergel bimodules \cite{Soe90}.

Polo proved that every polynomial with non-negative integer coefficients and constant term 1 appears as a Kazhdan--Lusztig polynomial of some sufficiently large
symmetric group \cite{Polo}.  
\end{example}

\begin{example}\label{ex:polytopes}
Let $\sigma$ be a pointed polyhedral cone and $P$ its poset of faces, ranked by codimension and ordered by reverse inclusion.
For all $F\leq G\in P$, let $\kappa_{FG}(t) = (t-1)^{r_{FG}}$; this defines a $P$-kernel precisely because the poset $P$ is Eulerian.
If $\Delta$ is a polytope obtained as a cross section of $\sigma$, then the polynomial 
$g_\Delta(t) := f_{\sigma 0}(t)$ is known as the {\bf \boldmath{$g$}-polynomial} of $\Delta$.
Since every interval in $P$ is isomorphic to the poset of faces of a smaller cone, all of the KLS-polynomials of $\sigma$
are $g$-polynomials of smaller polytopes.  (This contrasts with Example \ref{ex:coxeter}, where intervals in the Bruhat
poset need not be isomorphic to Bruhat posets for other groups.)

Suppose that $\sigma$ is rational, and let $X(\sigma)$ be the associated toric variety.  
This affine variety is stratified by torus orbits $O_F$ indexed by the faces of $\sigma$, with $\dim O_F = \rk(F)$ and 
$O_F\subset\overline{O}_G$ if and only if $F\leq G$.
For any point $p\in O_F$, Denef and Loeser \cite{DL} and Fieseler \cite{Fieseler} proved the following analogue of Equation \eqref{Weyl-IH}:
\begin{equation}\label{toric-IH}f_{FG}(t) = \sum_{i\geq 0} t^i \dim \IH^{2i}(\overline{O_G})_p.\end{equation}
In particular, taking $G=\{0\}$ and $F=\sigma$ so that $\overline{O_G} = X(\sigma)$ and $O_F = \{p\}$ is the unique fixed point, we have
$$g_\Delta(t) = \sum_{i\geq 0} t^i \dim \IH^{2i}(X(\sigma))_{p} = \sum_{i\geq 0} t^i \dim \IH^{2i}(X(\sigma)).\footnote{The
second equality follows from the fact that $X(\sigma)$ is an affine cone with cone point $p$.}$$
Non-negativity of the $g$-polynomial for arbitrary (not necessarily rational) polytopes was proved by Karu \cite{Karu}, using 
the theory of intersection cohomology of fans developed in \cite{BBFK}.

A full characterization of what polynomials occur as $g$-polynomials of polytopes is currently out of reach,
but the sequence of coefficients is conjecturally an $M$-sequence\footnote{An $M$-sequence is the sequence of graded dimensions of a commutative algebra that is generated in degree 1.} \cite{Kal:aspects,Stanley-IC}.
Evidence for this conjecture includes Bayer's generalized upper bound theorem \cite{Bayer}, which bounds the growth of the coefficients,
and Kalai's proof that the sequence of coefficients of $g_\Delta(t)$ has no internal zeros \cite[Theorem 1.4]{Braden-CICF}.
\end{example}

\begin{example}\label{ex:matroids}
Let $M$ be a loopless matroid on the ground set $E$, and let $\cL$ be the lattice of flats of $M$.
We can define an $\cL$-kernel by taking $\chi_{FG}(t)$ to be the characteristic polynomial of the interval $[F,G]\subset\cL$.
Rather than giving a definition of the characteristic polynomial, we note that the defining recursion for the KLS-polynomials associated with this kernel is equivalent
to the statement that the polynomials
\begin{equation}\label{eq:Z}Z_{F\!H}(t) := \sum_{F\leq G\leq H} t^{r_{F\!G}} f_{G\!H}(t)\end{equation} 
are ``palindromic'', that is, that $t^{r_{F\!H}}Z_{FH}(t^{-1}) = Z_{F\!H}(t)$.
This condition is relatively easy to understand, and it will also be the one that we use in Section \ref{sec:matroids}.
The polynomials $P_\cL(t) := f_{\emptyset E}(t)$ and $Z_\cL(t) := Z_{\emptyset E}(t)$ are called the {\bf Kazhdan--Lusztig polynomial} and
{\bf \boldmath{$Z$}-polynomial} of $\cL$ or $M$.
Since every interval in $\cL$ is itself the lattice of flats of a matroid, all of the KLS-polynomials of $\cL$
are Kazhdan--Lusztig polynomials of matroids.  In particular, for any $F\in\cL$, we denote the upper interval $[F,E]\subset \cL$ by $\cL_F$,
and Equation \eqref{eq:Z} translates to the identity
\begin{equation}\label{eq:Z2}Z_{\cL}(t) := \sum_{F\in\cL} t^{\rk F} P_{\cL_F}(t).\end{equation} 
See \cite{KLS} for a discussion of the analogues of $Z$-polynomials in the context of Examples \ref{ex:coxeter}
or \ref{ex:polytopes}.

Suppose that $M$ is the matroid associated with a linear subspace $L\subset\C^E$.  
Explicitly, this means that a subset $F\subset E$ is a flat if and only if there exists an element $p\in L$ such that $F = \{e\in E\mid p_e=0\}$.
The {\bf arrangement Schubert variety}\footnote{The word ``arrangement'' has to do with the fact
that $Y(L)$ is determined by the hyperplane arrangement in $L$ obtained by intersecting $L$ with each of the coordinate hyperplanes.} 
$Y(L)$ is defined to be the closure of $L$ inside of $(\C\P^1)^E$.  
The additive action of $L$ on itself extends to $Y(L)$, and the orbits $\{V_F\mid F\in\cL\}$ are indexed by flats.
We have $\dim V_F = \rk F$, and 
$V_F\subset\overline{V}_G$ if and only if $F\leq G$.  For $p\in V_F$, Elias, Proudfoot, and Wakefield \cite{EPW}
proved the analogue of Equations \eqref{Weyl-IH} and \eqref{toric-IH}:
\begin{equation}\label{matroid-IH}f_{FG}(t) = \sum_{i\geq 0} t^i \dim \IH^{2i}(\overline{V_G})_p.\end{equation}
In particular, taking $G=E$ and $F = \emptyset$ so that $V_\emptyset = \{\infty\}$, we have
$$P_\cL(t) = \sum_{i\geq 0} t^i \dim \IH^{2i}(Y(L))_{\infty}.$$
The $Z$-polynomial has a similar interpretation in terms of global intersection cohomology \cite{PXY}:
$$Z_\cL(t) = \sum_{i\geq 0} t^i \dim \IH^{2i}(Y(L)).$$
One can give analogous cohomological interpretations for matroids corresponding to linear spaces in positive characteristic, 
provided that one works with $l$-adic \'etale cohomology.
Non-negativity of these polynomials for arbitrary geometric lattices was proved by Braden, Huh, Matherne, Proudfoot, and Wang \cite{SHCG} by developing a theory of intersection
cohomology of matroids.  It was also reproved by Coron \cite{Coron} using different methods.  

Kazhdan--Lusztig polynomials of matroids are conjecturally real rooted \cite{kl-survey}; significant evidence
for this conjecture appears in \cite{fan-wheel-whirl, Chinese-uniform, stressed}.  This stands in 
stark contrast with Polo's theorem in Example \ref{ex:coxeter}.
\end{example}

\section{From topology to sheaves on posets.}\label{sec:topology}
There are two basic approaches to computing the intersection cohomology $\IH(X)$ of an algebraic variety $X$.  
First, one can follow the original definition of the intersection complex from \cite{GM2}, choosing a stratification of $X$ and
and adding one stratum at a time in order of increasing codimension.  Alternatively, one can find a resolution of singularities $\widetilde{X} \to X$; the decomposition theorem of \cite{BBD} implies that $\IH(X)$ is a direct summand of 
$\H(\widetilde{X})$ as a module over $\H(X)$.  In some cases, including Schubert and arrangement Schubert varieties, an argument of Ginzburg \cite{Gin-PS} shows that $\IH(X)$ is an indecomposable module, so the problem becomes finding a distinguished indecomposable summand of $\H(\widetilde{X})$.  

We will follow the first approach in this article, but the approach via the decomposition theorem has also played an important role in all three settings we discuss.  Indeed, the original definitions of Soergel bimodules \cite{Soe90,Soe92} and of the intersection cohomology module of a matroid \cite{SHCG} are as indecomposable direct summands of the cohomology of a resolution, and (a combinatorial version of) the decomposition theorem applied to a toric resolution of singularities plays a crucial role in Karu's proof \cite{Karu}.

\subsection{Torus equivariance and localization.}
It is often helpful to work equivariantly with respect to the action of a torus $T$.  For a toric variety $Y(\Delta)$, 
we take $T$ to be the natural torus that induces the orbit stratification;
 for a Schubert variety in the flag variety $G/B$, we take $T$ to be the maximal torus in the Borel subgroup $B\subset G$; 
 and for an arrangement Schubert variety $Y(L)$, we take $T=\cs$ acting on $L$ by homotheties.

The equivariant cohomology of a point is a polynomial ring
\[R = \H_T(\bullet) = \Sym(X^*(T)\otimes \R),\]
with the generating classes in degree $2$,
and the equivariant cohomology $\H_T(X)$ and equivariant intersection cohomology $\IH_T(X)$ are naturally $R$-modules.  The equivariant cohomology of a torus orbit $\cO$ is the quotient polynomial ring $$R_\cO := \H_T(\cO)\cong \H_{T_\cO}(\bullet),$$ where $T_\cO$ is the stabilizer of any point of $\cO$.  This ring is concentrated in even degrees, unlike the ordinary cohomology $H(\cO)$.  The action of $R$ on $\IH_T(X)_\cO$ factors through the quotient $R \to R_\cO$.

One advantage of working with equivariant intersection cohomology is that classes localize, meaning that they are determined by their restriction to local modules $\IH_T(X)_\cO$ at a finite number of $T$-orbits $\cO$.  Which $T$-orbits $\cO$ one uses depends on the type of variety one is considering.  
When $X$ is a Schubert variety or an arrangement Schubert variety (or an open union of strata in such a variety), 
restricting to the $T$-fixed points gives an injection
\[\IH_T(X) \hookrightarrow \bigoplus_{p\in X^T} \IH_T(X)_p.\]
This follows from the fact that
$\IH(X)$ vanishes in odd degrees, which also implies that $\IH_T(X)$ is {\bf equivariantly formal}: it is free as an $R$-module, and killing the equivariant parameters recovers the ordinary intersection cohomology $\IH(X)$.
For toric varieties, equivariant intersection cohomology classes also localize, but instead one restricts to every $T$-orbit:
\[\IH_T(X(\Sigma)) \hookrightarrow \bigoplus_{\sigma \in \Sigma} \IH_T(X(\Sigma))_{O_\sigma}.\]
(If $\Sigma$ is complete, or if it consists of a single full-dimensional cone and its faces, then $\IH_T(X(\Sigma))$ is equivariantly formal
and equivariant intersection cohomology localizes to the fixed points.)

In each case, the orbits being localized to are in bijection with the strata of $X$.  For (arrangement) Schubert varieties, each stratum is an affine space containing a unique $T$-fixed point, while for toric varieties the strata are the $T$-orbits themselves.  So computing $\IH_T(X)$ becomes a problem of computing the local modules $\IH_T(X)_\cO$ inductively, starting with an orbit in a smooth stratum and proceeding to more singular points.  

\subsection{Equivariant intersection cohomology as a sheaf.}
Consider a $T$-variety $X$ endowed with a $T$-invariant stratification $\{S_x\mid x\in P\}$ indexed by a poset $P$, where the order on $P$ is given by $a \le b$ if and only if $S_a \subset \overline{S_b}$.  
We say that $Q \subset P$ is an {\bf upper set} if, whenever $x\in Q$ and $x\leq y$, we also have $y\in Q$.
In this case,
\[U_Q:= \bigcup_{x\in Q} S_x\]
is an open subset of $X$, and all open unions of strata arise this way.  This means that we have a quotient map $q\colon X \to P$, where $P$ is made into a finite topological space with upper sets as open subsets.

Intersection cohomology has functorial restriction homomorphisms for open inclusions, so for open subsets $Q \subset Q'\subset P$, there is a homomorphism $\IH_T(U_{Q'})\to \IH_T(U_Q)$.  In other words, the assignment $$Q \mapsto \cF(Q) := \IH_T(U_Q)$$ defines a presheaf $\cF$ on $P$.  This can be viewed as the naive, non-derived pushforward of intersection cohomology along the quotient map $q$.  In general situations, one cannot hope for such a pushforward to be well-behaved; instead, one would need to work in the derived category and use spectral sequences and other tools from homological algebra for computations.  However, our varieties have two miraculous properties which help make this naive construction effectively computable.

The first miracle is that $\cF$ is actually a sheaf: for any cover $\{U_i\}$ of $U$
by open unions of strata, a collection of elements of $\IH_T(U_i)$ agreeing on the intersections will glue uniquely to an element of $\IH_T(U)$.  It's enough to do this for two sets $U = U_1 \cup U_2$, since our covers are always finite, and then this statement is equivalent to saying that the connecting homomorphisms in the Mayer-Vietoris sequence
\[ \dots \to \IH_T(U_1\cup U_2) \to \IH_T(U_1) \oplus \IH_T(U_2) \to \IH_T(U_1\cap U_2) \to \dots \]
vanish, so that the sequence becomes a short exact sequence of graded $R$-modules.

The second miracle is that this sheaf is \textbf{flabby}, meaning that the restriction map
$\IH_T(U_{Q})\to \IH_T(U_{Q'})$ is a surjection for any inclusion $Q' \subset Q$ of upper sets.
This is equivalent to the vanishing of connecting homomorphisms in the exact sequence
\[\dots \to \IH_T(U_Q,U_{Q'}) \to \IH_T(U_Q) \to \IH_T(U_{Q'}) \to \dots \]
of the pair $(U_Q,U_{Q'})$.  Flabbiness is important because it implies that the $R$-module $\cF(P)$ of global sections has a basis with nice triangularity properties: if $Q = Q' \cup \{x\}$, then a basis for $\IH_T(U_{Q'})$ lifts to a partial basis of $\IH_T(U_Q)$, with the new basis elements coming from $\IH_T(U_Q,U_Q')$.  Applying this recursively for all elements of $P$ yields a basis of $\cF(P)$ whose elements corresponding to a stratum $S_x$ will vanish on the complement of $\overline{S_x}$.

In both of the long exact sequences above, there are several ways to prove that the connecting homomorphisms vanish.  One way is to show that the mixed Hodge structures on the groups are pure of weight equal to the degree.  Alternatively, 
one can show that the intersection cohomology groups in question vanish in odd degrees.
In particular, using the decomposition theorem one finds that the following two geometric conditions are sufficient:
\begin{enumerate}
	\item[(1)] For all $x\in P$, pullback along $S_x \to \bullet$ gives a surjection 	$\H_T(\bullet) \to \H_T(S_x)$ in equivariant cohomology.
	\item[(2)] The variety $X$ admits an equivariant resolution whose fibers have only even cohomology.
\end{enumerate}
All of the varieties we consider satisfy both of these conditions.  Note that condition (1) would fail for toric varieties if we did not work equivariantly.

%

\subsection{The local computation.}
If $\cF$ is a sheaf on $P$, any section on an open set $Q\subset P$ is determined by its restrictions to the minimal open sets
\[P_x := \{y \in P \mid x \le y\}.\]
That means that, to construct $\cF$ inductively, it is enough to show how to extend $\cF$ to $P_x$ once it has been defined on $P^\circ_x := P_x \setminus \{x\}$.
In terms of open sets on the variety $X$, we set $U := U_{P_x}$ and $U^\circ := U \setminus S_x$, and we want to use our knowledge of $\IH_T(U^\circ)$ to compute $\IH_T(U)$.

A fundamental result of Bernstein and Lunts \cite[Section 14]{BL-equivariant} makes this possible.  
For simplicity, we assume that the new stratum $S_x = \{o\}$ is a single point fixed by $T$,
and that there is an affine neighborhood $N$ of $o$ and a one-dimensional subtorus $\C^* \subset T$ which contracts $N$ onto $\{o\}$;
see Remark \ref{slice} for a discussion of the more general case.
Let $N^\circ := N \setminus \{o\}$.

\begin{proposition}\label{prop:local calculation}
	The connecting homomorphism of the long exact sequence of $(N,N^\circ)$ in equivariant intersection cohomology vanishes, giving a short exact sequence
	\begin{equation}\label{eqn:local exact sequence}
		0\to \IH_T(N,N^\circ) \to \IH_T(N) \to \IH_T(N^\circ)\to 0.
	\end{equation}
	If $\dim X >0$, the following three statements hold: 
	\begin{enumerate}
		\item[{\em (a)}] $\IH_T(N)$ is a free $R$-module generated in degrees $0 \le d < \dim X$,
		\item[{\em (b)}] $\IH_T(N,N^\circ)$ is a free $R$-module generated in degrees $\dim X < d \le 2\dim X$, and 
		\item[{\em (c)}] $\IH_T(N)$ is canonically isomorphic to the minimal free $R$-module surjecting onto $\IH_T(N^\circ)$.
	\end{enumerate}
\end{proposition} 
In this proposition, the exact sequence \eqref{eqn:local exact sequence} and the statement (a) together imply the other statements.  Indeed, Poincar\'e--Verdier duality for intersection cohomology implies that (a) and (b) are equivalent, and 
(a) and (b) together with right exactness of the functor of reduction of scalars
\[M \mapsto \overline{M} := M \otimes_R \R\]
imply that $\overline{\IH_T(N)}\to \overline{\IH_T(N^\circ)}$ is an isomorphism.
That means that an $R$-basis of $\IH_T(N)$ maps to a minimal set of generators of $\IH_T(N^\circ)$.
The isomorphism in (c) is canonical because by (a) and (b) the generators of $\IH_T(N,N^\circ)$ are all in degrees higher than the generators of $\IH_T(N)$, so the only automorphism of $\IH_T(N)$ that commutes with the map to $\IH_T(N^\circ)$ is the identity.

Proposition \ref{prop:local calculation} allows us to compute $\IH_T(N)$ from $\IH_T(N^\circ)$.  This can be ``glued'' to $\IH_T(U^\circ)$ using the Mayer-Vietoris sequence for the cover $U = U^\circ \cup N$, which splits to give a short exact sequence
\[0\to \IH_T(U) \to \IH_T(U^\circ) \oplus \IH_T(N) \to \IH_T(N^\circ) \to 0.\]
Thus elements of $\IH_T(U)$ are given by pairs consisting of an element in $\IH_T(U^\circ)$ and an element of $\IH_T(N)$ that 
map to the same element in $\IH_T(N^\circ)$.  So we get an inductive computation of intersection cohomology, provided that one knows how to 
compute $\IH_T(N^\circ)$ from $\IH_T(U^\circ)$.  This problem has an answer in each of the cases we are considering, but the answers are particular to the geometry of each variety.  It would be very interesting to have a unified approach, but at present we do not.

\begin{remark}\label{slice}
We assumed above that the stratum $S_x$ is a point; here we briefly outline how to relax this assumption.
Assume instead that there is a point $p\in S_x$ such that the inclusion of the orbit $T\cdot p$ into $S_x$ is a homotopy equivalence.  (When $X$ is a toric variety, $S_x$ is itself a $T$-orbit, so $p$ can be any point.  When $X$ is a Schubert variety or an arrangement Schubert variety, $S_x$ is contractible, and 
we can take $p$ to be the unique $T$-fixed point in $S_x$.)  Then we define
\[R_x := \H_T(S_x) \cong \H_T(T\cdot p) = \H_{T_p}(\bullet),\]
where $T_p$ is the stabilizer of $p$.
We also assume that there exists a $T_p$-equivariant affine normal slice $N$ to $S_x$ at the point $p$, and that the stratum $S_x$ 
has a tubular neighborhood $\widehat{N}$ that is an $N$-bundle over $S_x$.  Letting $\widehat{N}^\circ = \widehat{N} \setminus S_x$, we have
\[\IH_T(\widehat{N}) \cong \IH_{T_p}(N),\]
so $\IH_T(\widehat{N})$ is the minimal free $R_x$-module surjecting onto $\IH_{T}(\widehat{N}^\circ)\cong \IH_{T_p}(N^\circ)$, and we proceed as before using the Mayer-Vietoris sequence of
$U = U^\circ \cup \widehat{N}$.
\end{remark}

\begin{remark}
Though we have chosen to focus on three classes of examples, we note that variations on this approach to computing intersection cohomology
have been employed in other situations, including hypertoric varieties \cite{TP08} and closures of $B\times B$-orbits in wonderful compactifications of adjoint groups \cite{Oloo}.
\end{remark}

\subsection{From spaces to posets.}\label{sec:no space}
We now rephrase the procedure above in the language of sheaves on posets.
Recall that we have defined the sheaf $\cF$ on $P$ by putting $\cF(Q) = \IH_T(U_Q)$ for any upper set $Q\subset P$.
Given an element $x\in P$, suppose that we have already computed $\cF(P_x^\circ) = \IH_T(U^\circ)$, and we want to 
compute $\cF(P_x) = \IH_T(U)$.
We first need to compute a ``boundary module'' $M_{\partial x} = M_{\partial x}(\cF)=\IH_{T}(\widehat{N}^\circ)$, which is a module over the polynomial ring $R_x$
defined in Remark \ref{slice}.  This is the step for which we have no uniform procedure,
but rather a recipe in each of the three families of examples that is specific to that family.
Once we have computed the $R_x$-module $M_{\partial x}$, take $M_x := \IH_T(N)$ to be the minimal free cover of $M_{\partial x}$
and we have a canonical isomorphism \[\cF(\Sigma_x) \cong \ker\Big(\cF(\Sigma^\circ_x)\oplus M_x \to M_{\partial x}\Big).\]
The induction begins when $x$ is a maximal element of $P$, so $U = S_x$ is a single open stratum.  
Since $U$ is smooth, the hypotheses of Remark \ref{slice} imply that
$\IH_T(U) = \H_T(U) = R_x$.  

The crucial point of this survey is that the resulting algorithm can be taken as a definition rather than a computation in the situations 
where the variety does not exist: nonrational polytopes/fans, 
Coxeter groups that are not Weyl groups, and non-realizable matroids.  But in each case it is very non-trivial to prove that the 
graded vector spaces that one constructs in this manner have Poincar\'e polynomials equal to the corresponding KLS-polynomials.
What is missing is an analogue of Proposition \ref{prop:local calculation}.  Specifically, one must show that $M_{\partial x}$ is generated in degrees less than 
$r_{xy}$, where $y$ is a maximal element of $P$.  When there really is a variety, this is a consequence of the hard Lefschetz theorem for intersection cohomology.  The existing proofs showing that combinatorially defined intersection cohomology has the expected degree vanishing properties all involve proving combinatorial analogues of the hard Lefschetz theorem and other results from Hodge theory, using complicated inductive schema.  But note that a recent result of Amini, Huh and Larson \cite{AHL} gives a more conceptual approach which avoids these complicated inductions, at least for matroids and polytopes.

\section{Polytopes and Fans.}\label{sec:fans}
In Sections \ref{sec:intro} and \ref{sec:KLS}, we proceeded in chronological order, beginning with Coxeter groups, moving on to 
polytopes (or cones or fans), and concluding with matroids.  Now, we will swap the order of the first two, treating polytopes before Coxeter groups.
This is partly because the theory of intersection cohomology of fans was the first to be fully worked out, and partly because it is simpler than the other two
theories, and is therefore a friendlier place to start.

\subsection{Rational fans.}
Let $\Sigma$ be a rational fan in $\R^d$.  We order the cones in $\Sigma$ by {\em reverse} inclusion, since smaller cones correspond to larger orbits in the toric variety $X(\Sigma)$.  This means that open sets in $\Sigma$ are subfans, which are sets of cones closed under taking faces.  The minimal open set $\Sigma_\tau$ containing $\tau$ consists of $\tau$ and all its faces, and 
$\Sigma^\circ_\tau$ is the set of all proper faces of $\tau$.

For any cone $\tau$ in $\Sigma$, the corresponding stratum $S_\tau$ is a quotient of the torus $T$, and its equivariant cohomology is canonically identified with $R_\tau = \Sym (\R\tau)^*$, the space of real-valued 
polynomial functions on $\tau$.  
These rings, together with natural restriction maps, fit together to form a sheaf of algebras $\cA$ on the poset $\Sigma$, so that the sections $\cA(\Sigma')$ on a subfan $\Sigma'\subset \Sigma$ are the real-valued functions on the support of $\Sigma'$ that are polynomial on each cone.  The sheaf $\cF$ induced by equivariant intersection cohomology is a sheaf of modules over $\cA$.  
Note that the sections $\cA(\Sigma')$ are not in general isomorphic to the equivariant cohomology ring of $X(\Sigma')$, but rather the equivariant Chow ring (tensored with $\R$).  This agrees with equivariant cohomology when $\Sigma'$ is simplicial or 
consists of a single cone and its faces.

The orbit $O_\tau$ corresponding to a cone $\tau \in \Sigma$ has an affine neighborhood that is again toric: it is the union of $O_\rho$ for all faces $\rho$ of $\tau$.    This implies that the boundary module $M_{\partial \tau}$ has a simple form: it is the module $\cF(\Sigma^\circ_\tau)$ of sections on the union of proper faces of $\tau$.  The general recipe given above says that the module of sections 
$\cF(\Sigma_\tau)$ for a new cone $\tau$ is the minimal free $R_\tau$-module that surjects onto $\cF(\Sigma^\circ_\tau)$.

\subsection{Arbitrary fans.}
As described in Section \ref{sec:no space}, one can make sense of the sheaf $\cF$ even if the fan $\Sigma$ is not rational.  In this case, there is no space,
so the procedure that we described for computing $\cF$ in the rational setting now becomes  a definition of $\cF$ for general fans.
Karu \cite{Karu} proved that, if $\Delta \subset \R^{d-1}$ is a polytope and $\Sigma\subset\R^d$ is the fan consisting of 
the cone over $\Delta \times \{1\}$ and all of its faces, then the $g$-polynomial $g_{\Delta}(t)$ is equal to the Poincar\'e polynomial (in $t^{1/2}$) of 
\[\overline{\cF(\Sigma)} := \cF(\Sigma) \otimes_R \R.\]

For simplicity, we only give examples of rational fans below, but these examples still provide some insight into the nature of Karu's argument
in the general case.

\begin{example}
The sheaf $\cA$ is flabby if and only if $\Sigma$ is a simplicial fan, meaning that each cone $\tau$ is generated by $\dim \tau$ rays.  
In that case, the equivariant intersection cohomology sheaf $\cF$ coincides with the structure sheaf $\cA$ itself.
If $\Sigma$ is rational, this reflects the fact that the intersection cohomology of an orbifold coincides with its cohomology.
For example, if $\Sigma$ is the fan in $\R^2$ whose maximal cones are the four orthants, and $R = \R[x,y]$ is the polynomial ring associated with any of the maximal cones, then 
	$\cF(\Sigma) = \cA(\Sigma)$ is the free $R$-module with basis consisting of the functions $1$, $|x|$, $|y|$ and $|xy|$.  
\end{example}

\begin{example}\label{square}
	For a non-simplicial example, we can lift the fan from the previous example to $\R^3$ by lifting the one-dimensional cones to the rays through $(1,0,1)$, $(-1,0,1)$, $(0,1,1)$, $(0,-1,1)$, and lifting the two-dimensional cones accordingly.  The resulting fan is the boundary of a $3$-dimensional cone $\hat\sigma$.  Let $\widehat{\Sigma}$ be the fan consisting of $\hat\sigma$ and all of its faces, and let $\hat\cF$ be the equivariant intersection cohomology sheaf on $\widehat{\Sigma}$.
	
	Since $\widehat{\Sigma}^\circ = \widehat{\Sigma} \setminus \{\hat\sigma\}$ is simplicial, we have 
	\[\hat\cF|_{\widehat\Sigma^\circ} = \cA|_{\widehat\Sigma^\circ},\]
	so $M_{\partial\hat\sigma}$ is the space of conewise polynomial functions on the boundary
    of $\hat\sigma$.  By projecting down to $\R^2$, we see that this is isomorphic to the space of sections $\cF(\Sigma)$ from the previous example, but now it is a module over the larger polynomial ring $\hat R = \R[x,y,z]$.  The multiplication by the new parameter $z$ is multiplication by $|x| + |y|$.  It follows that $M_{\partial\hat\sigma}$ is a non-free $\hat R$-module with minimal generating set $\{1,|x|\}$, so the structure sheaf $\cA$ of $\widehat{\Sigma}$ is not flabby.  Instead, $M_\sigma$ is    
    a free $\hat R$-module with generators in degree $0$ and $2$.  This corresponds to the fact that the $g$-polynomial of the square is $g^{}_{\square}(t) = 1+t$.
\end{example}

\begin{remark}
Example \ref{square} exhibits the general problem that Karu solved to show that the sheaf $\cF$ has stalks that categorify $g$-polynomials even when the fan $\Sigma$ is not rational.  The boundary of a full-dimensional cone $\hat\tau \subset \R^{d+1}$ can be ``flattened'' to give a complete fan $\Sigma$ in $\R^d$, and the piecewise linear function $\ell$ that lifts $\Sigma$ to $\partial\hat\tau$ represents an ample class in $\cA^2(\Sigma)$.  Karu showed that the hard Lefschetz theorem holds for the action of $\ell$ on $\cF(\Sigma)$, which implies that 
	$\hat\cF(\partial\hat\tau)$ is generated in degrees $< d+1$.  
	
	Karu's proof is a complicated induction starting from the case of a simplicial fan $\Sigma$, where hard Lefschetz was proved by McMullen \cite{McM}.  Recently, Adiprasito \cite{adiprasito2019} proved a hard Lefschetz result in the much more general setting of simplicial spheres.  It would be interesting to know if this could be extended to a theory of combinatorial intersection cohomology for regular CW-spheres that are not boundaries of convex polytopes.
\end{remark}

\section{Coxeter groups.}\label{sec:Coxeter}
We now turn to Coxeter groups and the theory of Soergel bimodules.  As in the previous section, we begin in the special case of a Weyl group,
where there is geometry, and the fundamental problem is to compute (rather than to define) certain intersection cohomology groups.

\subsection{Weyl groups.}
Let $G$ be a reductive group and $B\subset G$ a Borel subgroup, let $X := G/B$ be the flag variety, and let $X = \bigsqcup X_w$ be the Schubert stratification
indexed by the Weyl group $W$, as in Example \ref{ex:coxeter}.  Each Schubert cell $X_w$ is an affine cell with a unique $T$-fixed point, and we will abuse notation by also denoting this fixed
point by $w$.  We will fix throughout this section an element $w \in W$, and note that the strata of the Schubert variety $Y:= \overline{X_w}$ are indexed by the Bruhat interval $P := [e,w]\subset W$.  
Equivariant intersection cohomology of $\overline{X_w}$ induces a sheaf $\cF$ on $P$, 
whose sections on an upper set $Q \subset P$ are $\cF(Q) = \IH_T(U_Q)$.

These sections can be studied by equivariant localization.  
For all $v\in P$, let $M_{v} := \IH_T(Y)_v$ be the local equivariant intersection cohomology at the fixed point $v$.
We have an injective localization map $$\rho_Q = \bigoplus_{v \in Q} \rho_v \colon \IH_T(U_Q) \hookrightarrow \bigoplus_{v\in Q} M_{v}.$$ 
A result of Goresky, Kottwitz and MacPherson \cite{GKM} implies that the image of $\rho_Q$ is determined by information collected from the one-dimensional torus orbits, which can be expressed in terms of the moment graph $\Gamma$ for the flag variety.  This is a graph whose vertices are the $T$-fixed points $X^T$, and whose edges are the one-dimensional $T$-orbits in $X$.  The closure of a one-dimensional orbit is a projective line $\mathbb{P}^1$, and the two fixed points on this line are the vertices joined by the edge.   Each edge $E$ is labeled with a ``direction'' $\alpha_E \in X^*(T)$ whose kernel is the stabilizer of any point of $E$, so that $\H_T(E) \cong R/\alpha_E R$, where $R = \H_T(\bullet) = \Sym(X^*(T)\otimes \R)$.

The moment graph can be defined for any proper normal $T$-variety.\footnote{Given this fact, one might wonder why we work with the moment
graph for the flag variety rather for the Schubert variety.  The answer is that it is often useful to consider all of the Schubert varieties at the same time,
especially when studying Soergel bimodules, which we discuss briefly at the end of this section.}
In the particular case of the flag variety, it is known as the \textbf{Bruhat graph}, 
and it can be explicitly described as follows:
\begin{itemize}
	\item The vertex set may be identified with the Weyl group $W$.
	\item Two vertices $u,v\in W$ are joined by an edge if and only if $v = tu$ for some reflection $t\in W$.
	\item The direction $\alpha_E$ of an edge $w \stackrel{E}\longleftrightarrow tv$ is equal to the positive root corresponding to $t$.
\end{itemize}
For instance, if $G = GL_n$, the vertices are the symmetric group $S_n$, and edges are 
of the form $w \stackrel{E}\longleftrightarrow \tau_{ij}w$, where $\tau_{ij}$ is the transposition of $1\le i < j \le n$. The direction $\alpha_E$ of this edge is $e_i - e_j$.

Suppose that a pair of vertices $u,v\in P$ are connected by an edge $E$, and let $M_{E}:= \IH_T(Y)_E$.
We have a natural map \[M_{v} = \IH_T(Y)_v \to \IH_T(Y)_E = M_{E}\]
obtained by ``thickening'' the fixed point $x$ to a $T$-invariant affine neighborhood containing $E$ and then restricting to $E$.  
We now regard the moment graph $\Gamma$ as a poset whose elements are vertices and edges, 
with the only non-trivial order relations given by $v \unlhd E$ 
when a vertex $v$ lies on an edge $E$.\footnote{This partial order is not related to the Bruhat order on $W$.}
The modules $\{M_{v}\}$ and $\{M_{E}\}$, along with the restriction maps defined above, comprise the stalks of a sheaf $\cM$ on $\Gamma$, 
supported on the induced subgraph $\Gamma_{\!P}\subset \Gamma$ with vertex set $P\subset W$.
For any upper set $Q\subset P$, \cite[Section 6.3]{GKM} says that 
\[\cF(Q) = \IH_T(U_Q) \cong \cM(\Gamma_Q),\] the module of sections of $\cM$ on $\Gamma_{\!Q}$.\footnote{More precisely, we take the sections of $\cM$
on the minimal open subset of $\Gamma$ containing $\Gamma_{\!Q}$.}  Concretely, $\IH_T(U_Q)$ is isomorphic to the space of tuples $(m_v) \in \bigoplus_{v \in Q} M_{v}$ such that, for each edge $E$ containing vertices $u,v \in Q$, the images of $m_u$ and $m_v$ in $M_E$ are equal.

We can now describe the recursive computation of the sheaf $\cM$ on $\Gamma$, and therefore of the equivariant intersection cohomology sheaf $\cF$ on $P$.  The recursion starts by observing that 
\begin{equation}\label{eqn:B-M base case}
	M_{v} = 0 \;\, \text{for all} \;\, v \not \le w\;\; \text{and} \;\, M_{w} = R.	
\end{equation}
Next, we compute the edge modules $M_E$ for any edge $E$ joining vertices $u\leq v$, where $M_v$ has already been computed, via the formula
\begin{equation}\label{eqn:B-M constructibility}
	M_{E} \cong M_{v} \otimes_R R_E = M_{v}/\alpha_E M_{v}.
\end{equation}
Given a vertex $v\leq w$ of $\Gamma$, let $\Gamma_{> v}$ denote the smallest open subset of $\Gamma$ containing all of the vertices $\{u\in W\mid v<u\}$.
That is, it contains all such vertices, and all edges with at least one endpoint in this vertex set.  (It is typically not a graph, as it will contain ``dangling'' edges.)
The boundary module $M_{\partial v}$ is then determined by the sections on $\Gamma_{>v}$ as follows:
\begin{equation}\label{eqn:B-M boundary module}
	M_{\partial v} \cong \Im\left( \cM(\Gamma_{>v}) \to \bigoplus_E M_E\right),
\end{equation}
where the sum is over edges $E$ in $\Gamma_{>v}$ which are adjacent to $v$.
Then the module $M_v$ at the vertex $v$ is isomorphic to the minimal free $R$-module covering $M_{\partial v}$ \cite{BrM01}.

\begin{remark}
The condition \eqref{eqn:B-M base case} reflects the fact that the stratum $X_v$ is not contained in $Y = \overline{X_w}$ if $v \not\le w$, while $X_w$ is contained in the smooth locus of $Y$. The condition \eqref{eqn:B-M constructibility} comes from constructibility, since an edge $E$ connecting $v$ to a vertex which is smaller in the Bruhat order will be contained in the stratum $X_v$.  The surjectivity of $M_v \to M_{\partial v}$ reflects the fact that the equivariant intersection cohomology sheaf $\cF$ is flabby, since a section $\cF(P^\circ_v) = \cM(\Gamma_{>v})$ produces an element of $M_{\partial v}$, which can then be lifted to $M_v$, giving a section in $\cF(P_v)$.
Finally, 
the form the boundary module $M_{\partial w}$ takes in \eqref{eqn:B-M boundary module} is equivalent to the fact that the restriction
\begin{equation}\label{eqn:global sections to stalk}
\cF(P) = \cM(\Gamma) \to M_v	
\end{equation}
is surjective for all $v$, since given an element of $M_v$, its image in $M_{\partial v}$ is the image of a section in $\cM(\Gamma_{>v})$, and the resulting section of $\cF(P_v)$ extends to all of $P$ by flabbiness. 
The surjectivity of \eqref{eqn:global sections to stalk} is a special case of the fact that equivariant intersection cohomology is ``perfect'' (in the sense of Morse theory) for orderings coming from generic linear functions $\lambda\colon V \to \R$ \cite{Kirwan}. 
\end{remark}

\subsection{Arbitrary Coxeter groups.} If $W$ is an arbitrary Coxeter group, there is no flag manifold and no Schubert variety.  
However, for each element $w\in W$, we may still
define the poset $P = [e,w]\subset W$ and the sheaf $\cF$ on $P$ as above, where now the recursive procedure is taken as a definition rather than as
a computation.  If $W$ is a Weyl group, then 
Equation \eqref{Weyl-IH} from Example \ref{ex:coxeter} tells us that the Kazhdan--Lusztig polynomial $f_{vw}(t)$ is equal to the Poincar\'e polynomial of 
$\overline{M_v} \cong \IH(\overline{X_w})_v$.  For arbitrary $W$, the same fact follows 
from Elias and Williamson's proof of Soergel's conjecture \cite{EW14},
together with the relation between the moment graph sheaves and Soergel bimodules, which we now outline.

There is a natural sheaf $\cR$ of rings on $\Gamma$ with the property that $\cM$ is a sheaf of $\cR$-modules, and this sheaf has global section ring
\begin{equation*}\label{eqn:H_T(G/B)}\cR(\Gamma) = \{(f_v)_{v\in W}  \in \operatorname{Fun}(W,R) \mid f_u \equiv f_v \!\!\mod \alpha_E \;\;\text{for all edges}\;\; u \stackrel{E}{\longleftrightarrow} v\}.\end{equation*}
We then define $B_w := \cF(P) \cong\cM(\Gamma)$, which is a graded module over $\cR(\Gamma)$.
In the case where $W$ is a Weyl group, localization defines an isomorphism $\H_T(X)\cong \cR(\Gamma)$
\cite{GKM}, and the module structure is compatible with the isomorphism $B_w\cong \IH_T(\overline{X_w})$ described above.
In general, the ring $\cR(\Gamma)$ can be regarded as a purely algebraic stand-in for the equivariant cohomology of the flag manifold,
and the module $B_w$ as a stand-in for the equivariant intersection cohomology of a Schubert variety.

The ring $\cR(\Gamma)$ contains two copies of the polynomial ring $R$: one given by constant functions $v \mapsto f$, and one given by functions $v \mapsto f^v$, where $f^v$ is the ``twist'' of $f$ by $v$.  This makes $B_w$ into an $R$-bimodule, with the left (respectively right) action given by the first
(respectively second) inclusion $R\subset\cR(\Gamma)$.   
Fiebig \cite{Fie-comb-Cox-cats} showed that the resulting bimodules are isomorphic to the bimodules defined by Soergel \cite{Soe90}, which are the indecomposable direct summands of {\bf Bott--Samelson bimodules}
\[BS(s_1,...,s_r) := B_{s_1} \otimes \dots \otimes B_{s_r},\]
where the $s_i\in W$ are simple reflections and the tensor products are taken over $R$.
(When $W$ is a Weyl group, Bott--Samelson bimodules are isomorphic to equivariant cohomology rings of Bott--Samelson varieties, which resolve Schubert varieties.)

Soergel's approach to defining these bimodules has a number of advantages compared to the approach via moment graphs.  Most important
is the fact that, from the moment graph definition, the left and right actions of $R$ have a somewhat different nature, which makes the monoidal structure 
on the category of Soergel bimodules less accessible.  
The Grothendieck ring of this monoidal category is isomorphic to the Hecke algebra $\mathcal{H}_q(W)$; this categorification, known as the Hecke category,
has numerous applications in representation theory, knot invariants, and more.  

One place where moment graphs have an advantage is in the combinatorial invariance conjecture.   This is a longstanding conjecture saying that 
the Kazhdan--Lusztig polynomial $f_{vw}(t)$ only depends on the interval $[v,w]$ as an abstract poset.  Using moment graph sheaves, one 
can prove the weaker result that $f_{vw}(t)$ only depends on the induced subgraph $\Gamma_{[v,w]}$ along with the edge labelings.
A result of Dyer \cite{Dyer-Bruhat} implies that the graph $\Gamma_{[v,w]}$ is indeed determined by the isomorphism class of
the poset $[v,w]$, but the edge labelings are not.  The original conjecture remains open, but there has been renewed interest and activity around it recently, 
including \cite{Barkley-Gaetz,Williamson-comb-inv,Brenti-Marietti,Patimo}.

\subsection{Generalizations.} In our discussion so far we have confined ourselves to the case of the flag variety $G/B$ of a semisimple group $G$, but with only a few modifications the same story holds when $X$ is a partial flag variety $G/P$ of a Kac--Moody group $G$.  This includes many important applications in geometric representation theory which use the affine flag variety or affine Grassmannian, notably the geometric Satake theorem \cite{MV-GS}.
In this more general setting, the poset indexing the strata is in bijection with the quotient $W/W_P$ of the Weyl group by a parabolic subgroup.  This poset can be infinite, and $G/P$ can be infinite-dimensional, but lower intervals $[e,w]$ are still finite and Schubert varieties $\overline{X_w}$ are still finite-dimensional.    

Passing to Kac--Moody groups entails one significant complication, however.  The natural torus that acts may not be large enough to have a discrete set of one-dimensional orbits.  This means that the ``edges'' of the moment graph can connect more than two vertices, and the conditions they impose on equivariant cohomology and intersection cohomology become more complicated.  In terms of the moment graph, this is reflected in the fact that pairs of edges emerging from the same vertex can have parallel direction vectors, and in the Soergel bimodule setting it is reflected in the fact that the so-called ``geometric'' representation of $W$ coming from the root datum may not be ``reflection faithful'': some elements of $W$ that are not reflections can still act as reflections.  This causes problems both for moment graph sheaves and for Soergel bimodules.  

Soergel showed that there is a possibly larger representation of $W$ that is reflection faithful, so everything works as before; for the affine flag variety and affine Grassmannian it is enough to add the $\cs$ that acts by ``loop rotation''.  Later, Libedinsky \cite{Libedinsky} showed that the main theorems about Soergel bimodules (in particular the fact that their characters give Kazhdan--Lusztig polynomials) hold for the geometric representation if and only if they hold for Soergel's representation.  But this workaround is not completely satisfying, particularly when one tries to generalize from $\R$ coefficients to coefficients in a field of positive characteristic, where Soergel's trick does not work.  (We discuss the positive characteristic situation further in Section \ref{sec:positive} below.) One solution to this problem was given by Elias and Williamson \cite{EW-Soergel-calc}, who gave a diagrammatic description of the category of Soergel bimodules which is valid in a much broader setting.  By taking the idempotent completion one gets the desired Hecke category, although its objects are no longer bimodules.  An alternative approach was given by Abe \cite{Abe}, who constructs Soergel's category using $R$-bimodules equipped with localization maps to each element of $W$.

\section{Matroids.}\label{sec:matroids}
The last family of KLS-polynomials that we discuss is those coming from matroids, as in Example \ref{ex:matroids}.
Once again, we begin with the case of realizable matroids, in which there is a geometric object for us to study.

\subsection{Realizable matroids.}
Let $E$ be a finite set and $L\subset \C^E$ a linear subspace that is not contained in any coordinate hyperplane.  
The {\bf arrangement Schubert variety}
 of $Y(L)$, introduced in \cite{ArBoo}, 
is defined to be the closure of $L$ inside of $(\C\P^1)^E$.  Arrangement Schubert varieties are in many ways analogous to classical Schubert varieties, the most important shared feature being the existence of a stratification by affine varieties.  Indeed, for any subset $F\subset E$,
let $$V_F := \{p\in Y(L)\mid p_e \neq \infty\Leftrightarrow e\in F\}.$$
This locus turns out to be non-empty if and only if $F$ is a flat, and the dimension of $V_F$ is equal to the rank of $F$.

The action of $T=\cs$ on $L$ by scalar multiplication extends to an action on $Y(L)$.  Each stratum $V_F$ contains a unique fixed point $p_F$
whose coordinates are all equal to 0 or $\infty$, and $V_F\subset Y(L)$ is precisely the attracting set for $p_F$.
We also have $V_F = L\cdot p_F$, where the action of $L$ on $Y(L)$ is characterized by the property that it extends the additive action of $L$ on itself.
The stabilizer of $p_F$ in $L$ is equal to the linear subspace $L_F := L\cap \C^{E\setminus F}$, which implies that $V_F\cong L/L_F$ is isomorphic
to a vector space of dimension $\rk F$.
The closure of the stratum $V_G$ is the union of all $V_F$ such that $F\subset G$.

For any flat $F$, let $Y(L_F) = Y(L)\cap \big((\C\mathbb{P}^1)^{E\setminus F}\times\{0\}^F\big)$ be its associated arrangement Schubert variety.
This subvariety intersects the stratum $V_F$ transversely at the point $p_F$, and in fact $V_F$ has a $T$-stable Zariski open neighborhood 
 isomorphic to $V_F\times Y(L_F)$.  Thus we may think of $Y(L_F)$ as a normal slice to the stratum $V_F$.
The interval $\cL_F = [F,E]\subset\cL$ is canonically identified with 
the poset  of flats for $L_F$.

Since the torus $T$ is only 1-dimensional, there are always infinitely many 1-dimensional $T$-orbits (unless $\dim L$ is equal to 0 or 1).  Thus $Y(L)$ does not induce a moment graph as in the previous section.
Nevertheless, there is a simple formula for its equivariant cohomology:  $\H_T(Y(L))$ is a free module over the ring $R = \H_T(\bullet) = \R[\hbar]$ with
basis $\{y_F\}_{F\in \cL}$ and multiplication rule
\[y_F y_G = \hbar^{\rk F + \rk G - \rk F\vee G} y_{F \vee G}.\]
Here the {\bf join} $F\vee G$ is equal to the smallest flat containing both $F$ and $G$.   
The basis element $y_{\emptyset}$ is equal to the multiplicative identity, and the degree of $y_F$ is equal to $2\rk F$.
The restriction homomorphism $$\rho_F:\H_T(Y(L))\to \H_T(p_F)\cong R$$ takes $y_G$ to $\hbar^{\rk G}$ if $G\subset F$ and to 0 otherwise.
The restriction $$\varphi_F:\H_T(Y(L))\to \H_T(Y(L_F))$$ takes $y_G$ to 
$\hbar^{\rk F + \rk G -\rk F\vee G}y_{(F\vee G)\setminus F}$.

Let $\cA$ be the sheaf of algebras on the poset $\cL$ whose ring of sections on any open subset of $\cL$ is the 
equivariant cohomology ring of the corresponding union of strata.  This sheaf is characterized by the properties that $\cA(\cL_F) = \H_T(Y(L_F))$
and the restriction map from global sections to sections on $\cL_F$ is given by $\varphi_F$.

Let $\cF$ be the sheaf of $\cA$-modules whose module of sections on any open subset of $\cL$ is the 
equivariant intersection cohomology of the corresponding union of strata.
This sheaf can be computed inductively as follows.  We have $M_E = \cF(\{E\}) = \cA(\{E\}) = R$.  If we have already computed $\cF$ on $\cL^\circ_F = \cL_F \setminus \{F\}$, we let
\[M_{\partial F} := \cF(\cL^\circ_F) \otimes_{\cA(\cL_F)} R,\]
where $\cA(\cL_F)$ acts on $\cF(\cL^\circ_F)$ by restriction and $\cA(\cL_F)$ maps to $R$ via $\rho_F$.  Equivalently, $M_{\partial F}$ is the quotient of $\cF(\cL^\circ_F)$ by the action of $y_G$ for $G \not\le F$.  If $M_F$ is the minimal free $R$-module surjecting onto $M_{\partial F}$, then
$\cF(\cL_F)$ is canonically isomorphic to the kernel of the map
\[\cF(\cL_F^\circ) \oplus M_F \to M_{\partial F}.\]
Repeating this step for all flats in order of decreasing rank computes the sheaf $\cF$ on the entire poset $\cL$.

\begin{example}
Suppose that $|E|=1$ and $L = \C$, so that $Y(L) = \C\P^1$.  The $R$-algebra $\H_T(Y(L))$ has an additive basis consisting of $y_\emptyset$ and $y_E$,
where $y_\emptyset$ is the multiplicative identity and $y_E^2 = \hbar y_E$.  We begin with the isomorphisms $\cL_\emptyset^\circ=\{E\}$ and $\cF(\{E\}) = \cA(\{E\}) = R$.
We have two different homomorphisms from $\H_T(Y(L))$ to $R$:  one is $\rho_\emptyset$, which 
takes $y_\emptyset$ to $1$ and $y_E$ to $0$; the other is the restriction map from $\cH_T(\cL)$ to $\cH_T(\{E\})$, which takes
$y_\emptyset$ to $1$ and $y_E$ to $\hbar$.
We define $$M_{\partial\emptyset} := R\otimes_{\H_T(Y(L))} R \cong \R,$$
where the tensor product uses the two different homomorphisms from $\H_T(Y(L))$ to $R$.
We then define $M_\emptyset$ to be the minimal free cover of $\R$, which is isomorphic to $R$,
and $\cF(\cL)$ is identified with the kernel of the graded
$R$-algebra homomorphism $R\oplus R\to\R$.  This is a free $R$-module with generators in degrees 0 and 1, and is isomorphic to the regular
$\cA(\cL) = \H_T(Y(L))$-module.  This is as it should be, since $\C\P^1$ is smooth, and intersection cohomology is therefore isomorphic to cohomology.
\end{example}

\begin{example}
Suppose that $|E|=4$ and $L\subset\C^E$ is the 3-dimensional subspace consisting of vectors whose coordinates add to zero.  Its flats are all subsets $F \subset E$ except
subsets with $|F| = 3$. 
This is the smallest example for which the arrangement Schubert variety has worse than orbifold singularities, and therefore the first example
for which intersection cohomology does not coincide with cohomology.  
For every $F\neq \emptyset$, we find that $M_F = R$ and $\cF(\cL_F) = \cA(\cL_F)$,
which reflects the fact that $Y(L)$ has only orbifold singularities away from the most singular point.
However, the calculation at the most singular point is more subtle.  The graded $R$-module $M_{\partial\emptyset}$ has the property
that $\overline{M_{\partial\emptyset}}\cong \R \oplus\R^2[-1]$, hence $M_{\partial\emptyset}$ 
has a minimal free cover $M_\emptyset\cong R\oplus R^2[-1]$.  This reflects the fact that the local intersection cohomology Betti numbers at the most singular
point are 1 and 2.
This allows us to compute $\cF(\cL) \cong R \oplus R^6[-1] \oplus R^6[-2]\oplus R[-3]$,
which shows us that the intersection cohomology Betti numbers of $Y(L)$ are 1, 6, 6, and 1.
\end{example}

\subsection{Arbitrary matroids.}\label{sec:arbitrary}
Now consider the case of an arbitrary loopless matroid with lattice of flats $\cL$.\footnote{Since the geometric lattice $\cL$ contains all of the information of the matroid up to simplification, there is no need for us to give our matroid a name.
This allows us to reserve the letter $M$ for the modules that arise in our sheaf theoretic constructions.}
We may define the sheaf $\cA$ of rings on $\cL$ purely algebraically, and construct the sheaf $\cF$ of $\cA$-modules as above.
If $\cL$ comes from a linear subspace $L\subset\C^E$, then we have canonical isomorphisms
$$M_F \cong \IH_T(Y(L_F))_{p_F}\qquad\text{and}\qquad \cF(\cL)\cong \IH_T(Y(L)),$$
and therefore 
$$\overline{M_F} \cong \IH(Y(L_F))_{p_F}\qquad\text{and}\qquad \overline{\cF(\cL)}\cong \IH(Y(L)).$$
By Example \ref{ex:matroids}, this implies that the matroid Kazhdan--Lusztig polynomial $P_{\cL_F}(t)$ is equal to the Poincar\'e polynomial of $\overline{M_F}$,
and $Z_{\cL}(t)$ is equal to the Poincar\'e polynomial of $\overline{\cF(\cL)}$.  But in fact they are equal for arbitrary matroids, whether or not they come from a linear subspace $L$.  To show this, we need to prove
four statements:
\begin{itemize}
\item[(1)] If $E$ is the maximal element of $\cL$, then $\overline{M_E}\cong\R$.
\item[(2)] For any $F\neq E$, the graded vector space $\overline{M_F}$ vanishes in degrees greater than or equal to $\rk E - \rk F$.\footnote{Note
that we are using the convention in which the degree of $y_F$ is equal to $2\rk F$.  If we used the convention that the degree of $y_F$ was equal to $\rk F$,
then $\overline{M_F}$ would have to vanish in degrees greater than or equal to $(\rk E - \rk F)/2$.}
\item[(3)] For any $k$, we have $\dim \overline{\cF(\cL)}^{2k} = \dim \overline{\cF(\cL)}^{2(\rk E - k)}$.
\item[(4)] There exists an isomorphism of graded vector spaces $\displaystyle\overline{\cF(\cL)} \cong \bigoplus_F \overline{M_F}[-2\rk F]$.
\end{itemize}
Statements (1) and (2) correspond to the first two conditions in the definition of KLS-polynomials from Section \ref{sec:KLS},
while statements (3) and (4) correspond to palindromicity of the Z-polynomial as defined in Equation \eqref{eq:Z2}.

Statement (1) is trivial, and statement (4) can be proved by choosing an $R$-basis for $M_F$ for every $F$ and 
lifting all of these basis elements to an $R$-basis for $\cF(\cL)$.
Statement (3) reflects the fact that the dual module $\cF(\cL)^*$ is isomorphic to $\cF(\cL)$ up to a shift.  
When the variety $Y(L)$ exists, this is a manifestation of Poincar\'e duality for intersection cohomology.  
In the general case, this can be proved by considering  triangular ``flow-up'' and ``flow-down'' bases for $\cF(\cL)$ with dual support properties.  


Statement (2) is the most subtle.  When $\cL$ is realized by a vector subspace $L$, it follows from the fundamental local calculation 
in Proposition \ref{prop:local calculation}, which is a consequence of the fact that the intersection cohomology of a projective variety satisfies 
the hard Lefschetz theorem.  
Just like the proofs in \cite{Karu} for nonrational fans and \cite{EW14} for Coxeter groups that are not Weyl groups, 
the proof of this property for general matroids in \cite{SHCG} involves a complicated induction simulating various Hodge-theoretic 
statements including the hard Lefschetz theorem and the Hodge--Riemann bilinear relations.

In the paper \cite{SHCG}, we follow a different but related approach to constructing the intersection cohomology of 
a matroid.  We work non-equivariantly, constructing a module $\IH(\cL)$ over the \textbf{graded M\"obius algebra} $\H(\cL) := \overline{\cA(\cL)}$.  
Then we find $\IH(\cL)$ 
as a direct summand of an $\H(\cL)$-algebra $\CH(\cL)$ known as the \textbf{augmented Chow ring} of $\cL$, 
which in the case of a subspace $L\subset\C^E$ may be identified with the cohomology of a resolution of $Y(L)$.
Thus the construction in \cite{SHCG} is closer to Soergel's definition of the bimodules $B_w$ in that it involves simulating the decomposition theorem rather than constructing intersection cohomology inductively stratum-by-stratum or flat-by-flat.

These approaches are related by the existence of a canonical isomorphism $\IH(\cL) \cong \overline{\cF(\cL)}$.  This isomorphism appears in the forthcoming papers \cite{BHMPW-IHT-module,BHMPW-IH-module}.  The paper \cite{BHMPW-IH-module} gives simple conditions which uniquely characterize $\IH(\cL)$, and \cite{BHMPW-IHT-module} shows that $\overline{\cF(\cL)}$ satisfies these conditions.  

\section{Inequalities.}
The most obvious application of constructing a graded vector space whose Poincar\'e polynomial is equal to a given KLS-polynomial is to show that the polynomial has non-negative coefficients.  However, it can also lead to more subtle inequalities, as well.  

\subsection{Monotonicity.}
For polynomials $f(t),g(t)\in \Z[t]$, we write $f(t) \preceq g(t)$ if $g(t)-f(t)$ has nonnegative coefficients.  If $\{f_{xy}(t)\mid x\leq y \in P\}$ is one of the three families of KLS-polynomials we are considering, we have 
\begin{equation}\label{eqn:monotonicity}
	f_{xz}(t) \preceq f_{yz}(t) \;\; \text{whenever} \; x \le y \le z.
\end{equation}
This inequality was first proved for Kazhdan--Lusztig polynomials of Weyl groups by Irving \cite{Irv88}, using the interpretation of their coefficients as multiplicities of simple modules in Verma modules of the associated Lie algebra.  
A geometric proof for finite and affine Weyl groups appeared in \cite[Corollary 3.7]{BrM01}, 
and this proof generalizes to all three settings in which we have a variety.  We briefly summarize the argument here.

Let $\cF$ be the equivariant intersection cohomology sheaf on $P^z = \{x\in P\mid x\leq z\}$.  This means that, for each $x\leq z$, we have
a graded $R_x$-module $M_x$ with the property that the Poincar\'e polynomial of $\overline{M_x} = M_x\otimes_{R_x}\R$ 
is equal to the KLS-polynomial $f_{xz}(t)$.
For any $x \le z$, the restriction homomorphism \[\varphi_x\colon \cF(P^z) \to M_x\] 	
is surjective, since elements of $M_x$ can first be extended to $[x,z]$, and then extended to all of $P^z$ by flabbiness.  The reduced homomorphism
$\overline{\varphi_x}\colon \overline{\cF(P^z)} \to \overline{M_x}$ is therefore also surjective.

We claim that, for any $x \le y \le z$, there is a unique homomorphism $\psi_{yx}\colon \overline{M_x} \to \overline{M_y}$ such that $\overline{\varphi_y} = \psi_{yx}\overline{\varphi_x}$.	
Indeed, if $\cG \in D^b(X)$ is constructible with respect to a stratification $X = \bigsqcup_{x} S_x$ with connected strata and $S_x \subset \overline{S_y}$, then there is a homomorphism from stalk cohomology of $\cG$ at a point of $S_x$ to stalk cohomology at a point of $S_y$, commuting with the restriction maps from global cohomology $\H(\cG)$.
Since $\overline{\varphi_y}$ is surjective, $\psi_{yx}$ must be surjective too, which gives the inequality \eqref{eqn:monotonicity}.

The existence of the maps $\psi_{yx}$, and therefore the inequality \eqref{eqn:monotonicity}, has been proved for arbitrary Coxeter groups \cite{Plaza} 
and for arbitrary matroids \cite[Theorem 1.4]{SHCG}, both using machinery specific to their settings. 
The inequality \eqref{eqn:monotonicity} does hold for $g$-polynomials of polytopes, but in fact a stronger inequality holds: for any pointed cone $\sigma$ and any face $\tau$ of $\sigma$, we have
\begin{equation}\label{eqn:Kalai}
	f_{\sigma 0}(t) \succeq f_{\sigma\tau}(t)f_{\tau 0}(t).	
\end{equation}
This implies the inequality \eqref{eqn:monotonicity} because $f_{\tau 0}(t) \succeq 1$.  

The inequality \eqref{eqn:Kalai} was originally conjectured by Kalai \cite{Kal:new_basis}, who proved it for the coefficients of $t$ and $t^2$ \cite{Kal:new_basis,Kal:aspects}.  
A proof for rational cones via the geometry of toric varieties appeared in \cite{BrM:Kconj}.  
The key geometric fact that was used is the following: if $X(\sigma)$ is the toric variety defined by a rational cone $\sigma$, and $Y = \overline{S_\tau}$ is the closure of the orbit corresponding to a face $\tau$, then  there is a one-dimensional subtorus $T_0\subset T$ that contracts $X(\sigma)$ onto the fixed point locus $Y = X(\sigma)^{T_0}$.  This implies that the restriction to $Y$ of $\IC_{X(\sigma)}$ is isomorphic to a direct sum of shifted 
IC sheaves of orbit closures; the number of copies of $\IC_Y$ is equal to $f_{\sigma \tau}(t)$, with the degree measuring the shifts.
Taking cohomology of the stalks at the unique $T$-fixed point, $\IC_{X(\sigma)}|_Y$ contributes $f_{\sigma 0}(t)$, $\IC_{Y}$ contributes $f_{\tau 0}(t)$,
and we obtain the desired inequality.

It was noted in \cite{BBFK,BrLu} that the same argument works in the setting of combinatorial intersection cohomology of fans: 
the restriction of the sheaf $\cF_0$ to the interval $[\sigma, \tau]$ 
is flabby and has free stalks, which implies that it is isomorphic to a direct sum of shifts of sheaves $\cF_\nu$ for cones $\nu \le \tau$.  So Karu's theorem \cite{Karu} implies that \eqref{eqn:Kalai} holds for non-rational cones.  

\begin{remark}The torus actions on Schubert varieties and arrangement Schubert varieties do not have the analogous property; indeed the closure of a stratum may not be the fixed-point locus of a subtorus of $T$.  The inequality \eqref{eqn:Kalai} is false for classical Kazhdan--Lusztig polynomials.  
The analogous inequality for Kazhdan--Lusztig polynomials of matroids has been conjectured \cite[Conjecture 8.22]{EMPV},
but any proof will require new ideas. 
\end{remark}

\subsection{Top-heaviness.}
Another class of inequalities which can be deduced from these constructions only applies to the Coxeter and matroid settings.  We begin with the case where there is a variety.  Let $Y$ be either a Schubert variety $\overline{X_y}$ or an arrangement Schubert variety $Y(L)$ associated with a linear space $L \subset \C^E$, and let $P$ be the poset indexing strata of $Y$, that is, the interval $[e, y]\subset W$ or the lattice of flats $\cL$ of $L$.
Let 
\[	h_j := \#\{x \in P \mid \rk x = j\},\]
and let $d = \dim Y$ be the rank of the maximal element of $P$.
The poset $P$ is {\bf top-heavy} in the following sense:
\begin{equation}\label{eqn:top-heavy}
h_j \le h_{k}\;\;\text{whenever}\;\; j \le k \le d-j.
\end{equation}
In particular, $h_j \le h_{d-j}$ for $j \le d/2$, and the sequence is increasing for the first half: $1 = h_0 \le h_1 \le \dots \le h_{\lfloor d/2\rfloor}$. 

Bj\"orner and Ekedahl \cite{BE09} proved \eqref{eqn:top-heavy} for Schubert varieties.  Since the stratum $S_x$ is an affine space of dimension equal to the rank of $x$ for any $x\in P$, we have
$h_j = \dim \H^{2j}(Y)$.  Then an argument using mixed Hodge structures shows that the natural homomorphism $\H(Y) \to \IH(Y)$ is injective.  Since $Y$ is projective, intersection cohomology satisfies the hard Lefschetz theorem: there is a class $\a\in \H^2(Y)$ so that 
\[\a^{d-2j}\colon \IH^{2j}(Y) \to \IH^{2(d-j)}(Y)\]
is an isomorphism for $j \le d/2$, and so $\a^{k-j}\colon\IH^{2j}(Y) \to \IH^{2k}(Y)$
is an injection for $j \le k \le d-j$.
Since $\H(Y) \to \IH(Y)$ is a homomorphism of $\H(Y)$-modules, it follows that 
\[\a^{k-j}\colon\H^{2j}(Y) \to \H^{2k}(Y)\]
is injective, giving \eqref{eqn:top-heavy}.  
Huh and Wang \cite{HW} showed that the same argument applies to arrangement Schubert varieties, which gives top-heaviness for geometric lattices.

The case of general geometric lattices is proved in \cite{SHCG}, by showing that the graded M\"obius algebra $\H(\cL)$ injects into the combinatorial 
intersection cohomology module $\IH(\cL)$ and that $\IH(\cL)$ satisfies 
the hard Lefschetz property; this settles a 1975 conjecture of Dowling and Wilson \cite{DW75}.  The case of Coxeter groups that are not Weyl groups can be proved using Soergel bimodules, and appears in an unpublished preprint of Melvin and Slofstra.

\section{Positive characteristic.}\label{sec:positive}
Up to this point, we have assumed that the coefficients for our intersection cohomology groups and sheaves are taken in $\R$;
we have done this because crucial results from Hodge theory such as the decomposition theorem and the hard Lefschetz theorem require coefficients in a field of characteristic zero.  These results are responsible for the necessary parity vanishing properties of 
our various cohomology groups.

What happens when we try to work over a field of characteristic $p>0$?
Intersection cohomology can still be defined, but it is very difficult to compute; in particular,
it can be nonzero in odd degrees.
Furthermore, the direct summands of the derived pushforward of the constant sheaf from a resolution
may not be intersection complexes, even if the fibers of the resolution have no odd cohomology.
%

An alternative is to work with \textbf{parity sheaves}, as defined by 
Juteau, Mautner and Williamson \cite{JMW14}.  For a stratified variety $X$ satisfying appropriate hypotheses, each irreducible local system on a stratum $X_y$ has a unique indecomposable extension that is supported on $\overline{X_y}$, with stalks and costalks whose cohomology vanishes in odd degrees.  It satisfies a version of Poincar\'e--Verdier duality: up to a shift, the Verdier dual of a parity sheaf is isomorphic to the parity sheaf associated with the dual local system.  There is also a version of the decomposition theorem for parity sheaves: for a resolution of singularities whose fibers have no odd cohomology, the derived pushforward of the constant sheaf
decomposes as a direct sum of shifts of parity sheaves. 

This decomposition is determined by certain intersection forms on the homology groups of the fibers of the resolution.
For all but finitely many primes $p$, the ranks of these forms when taken over a field of characteristic $p$ 
are the same as the ranks over $\R$, the parity sheaves are intersection complexes, and the Poincar\'e polynomials
of their stalks coincide with KLS-polynomials.
However, for some special characteristics, the ranks of the forms can drop, and the stalk cohomology 
of the parity sheaves can be nonzero in high degree.  Thus the associated 
{\bf \boldmath{$p$}-KLS-polynomials} ${}^p\!f_{xy}(t)$, which are Poincar\'e polynomials of stalks of parity sheaves in analogy with Equations \eqref{Weyl-IH}, \eqref{toric-IH}, and \eqref{matroid-IH},
need not have degrees strictly less than $r_{xy}/2$.  Since this degree restriction is what makes the KLS recursion at the beginning
of Section \ref{sec:KLS} give uniquely defined polynomials, computing the $p$-KLS polynomials is much more difficult.

\begin{remark}
Parity sheaves on toric varieties have not been extensively studied, perhaps because the $p$-KLS-polynomials 
depend on the lattice which respect to which the fan is rational.  In particular, 
it doesn't make sense to ask for a positive characteristic theory for non-rational fans.
\end{remark}

\subsection{Flag varieties and \boldmath{$p$}-Kazhdan--Lusztig polynomials.}
For finite and affine flag varieties with the Schubert stratification, parity sheaves had already been considered by Soergel \cite{Soergel-modular} and Fiebig \cite{Fiebig-Lusztig-conj}, who showed that direct summands of derived pushforwards of constant sheaves from Bott--Samelson resolutions of Schubert varieties give information about representations of the Langlands dual group in positive characteristic.  Later it was shown that parity sheaves on the affine Grassmannian correspond under the geometric Satake equivalence to tilting modules for the Langlands dual group for all but a small number of characteristics $p$ \cite{JMW-tilting,MR}.  As a result, the $p$-Kazhdan--Lusztig polynomials and the $p$-canonical basis of Hecke algebras have been of great interest for recent developments in modular representation theory.

Moment graphs provide one way to compute $p$-Kazhdan--Lustig polynomials.  
The parity vanishing properties of parity sheaves imply the same vanishing of connecting homomorphisms that featured in the computatation of equivariant intersection cohomology outlined in Section \ref{sec:topology}, and the local computation provided by Proposition \ref{prop:local calculation} still works, except that the degree restrictions in (a) and (b) may not hold, and the isomorphism in (c) may not be canonical. 
Fiebig and Williamson \cite{FW} showed that, in positive characteristic, the moment graph algorithm outlined in Section \ref{sec:Coxeter} actually computes the cohomology of the $T$-equivariant parity sheaf $\cE_w$ supported on $\overline{X_w}$ (since the stratum $X_w$ is simply connected, there is only one irreducible local system on $X_w$, and so only one parity sheaf with this support).  

More precisely, this is true provided that the graph is ``$p$-GKM'', which means that the labels on edges adjacent to the same vertex do not become parallel modulo $p$.  For finite Weyl groups this rules out a few small primes, but for affine Weyl groups the restriction is more serious.  
The graph of a Schubert variety in the affine flag variety will be $p$-GKM for all but finitely many $p$, but the list of bad primes grows as the Schubert variety gets larger, and there are no characteristics that work for the graph of the entire affine flag variety. 

Soergel bimodules give another path to computing $p$-Kazhdan-Lusztig polynomials, which has proved more effective in practice.  
There is a technical difficulty arising from the fact that there is no reflection faithful representation of the affine Weyl group in positive characteristic.  But the diagrammatic calculation of homomorphisms between Bott--Samelson bimodules remains valid in any characteristic, which is enough to compute the necessary intersection forms and projection maps.  This makes $p$-Kazhdan--Lusztig polynomials effectively computable \cite{GJW,JW-p-canon}, although the computation remains much more difficult than that of the classical (characteristic zero) Kazdhan--Lusztig polynomials.

This circle of ideas has been successfully applied to important problems in modular representation theory, particularly to the representations of algebraic groups in positive characteristic.  Lusztig's conjecture on the characters of irreducible representations was previously known to hold for all sufficiently large characteristics depending on the group, but, using a combination of moment graph and bimodule techniques, Fiebig found an explicit (enormous) bound for the bad characteristics \cite{Fiebig-Lusztig-conj,Fiebig-bound}.  On the other hand, Williamson used a small subset of the diagrammatic calculus to show that the bad primes for $G = \operatorname{SL}_n$ become exponentially large as $n$ grows, disproving the conjecture that Lusztig's formula would hold for all primes above the Coxeter number.  More recently, Riche and Williamson \cite{RW-tilting} conjectured a character formula for tilting modules involving the $p$-Kazhdan--Lusztig polynomials. They proved this conjecture for $\operatorname{GL}_n$, and it was extended to all groups by Achar, Makisumi, Riche, and Williamson \cite{AMRW}.

\subsection{Matroids.}\label{p-matroids}
The story for arrangement Schubert varieties and matroids is very similar, but somewhat simpler. 
Let $L\subset \C^E$ be a linear subspace and $\cL$ the lattice of flats.
The arrangement Schubert variety $Y(L)$ with its orbit stratification satisfies the hypotheses for the existence of $T$-equivariant parity sheaves with coefficients in any field.  Since the strata are contractible, there is one parity sheaf (up to shift) for each stratum.  Since the closure of a stratum is isomorphic to another arrangement Schubert variety, it is enough to consider the parity sheaf $\cE$ supported on all of $Y$.  

As in the characteristic zero case, the parity vanishing conditions imply that $\cE$ induces a sheaf $\cF$ on the poset $\cL$, whose sections on an open subset of $\cL$ coincide with the equivariant cohomology of $\cE$ on the corresponding open union of strata.  The algorithm described in Section \ref{sec:matroids} still works to compute $\cF$.  The 
\textbf{\boldmath{$p$}-Kazhdan--Lusztig polynomial} ${}^p\!P_{\cL}(t)$ is defined to be the Poincar\'e polynomial of the reduced stalk module $\overline{M_\emptyset} = \H(\cE|_{p_\emptyset})$, and the \textbf{\boldmath{$p$}-\boldmath{$Z$}-polynomial} ${}^p\!Z_{\cL}(t)$ is defined to be the  
Poincar\'e polynomial of the reduced global sections $\overline{\cF(\cL)} = \H(\cE)$.
As in Section \ref{sec:matroids}, we can make sense of this definition for an arbitrary matroid, where the variety and the parity sheaf no longer exist,
but we still have a sheaf $\cF$ on the poset $\cL$.
The properties (1),(3), and (4) of Section \ref{sec:arbitrary} still hold, so these polynomials satisfy the relation
\[{}^p\!Z_{\cL}(t) = \sum_{F \in \cL} t^{\rk F}\,{}^p\!P_{\cL_F}(t),\]
and the polynomial ${}^p\!Z_{\cL}(t)$ is palindromic of degree equal to the rank of $\cL$.  
But ${}^p\!P_{\cL}(t)$ need not have degree less than half of the rank, so these properties no longer uniquely determine the polynomials.    

At present, the only known ways to compute ${}^p\!P_{\cL}(t)$ are to compute the sheaf $\cF$ (or a non-equivariant analogue),
or to compute certain intersection forms on stalks of the augmented Chow ring in characteristic $p$, which is a purely algebraic stand-in
for the cohomology ring of a resolution.
But there is evidence to suggest that the answer will be simpler than it is in the case of Coxeter groups.  
For instance, we have a simple criterion for when our polynomials are trivial.
The lattice $\cL$ is said to be \textbf{modular} if 
$\cL$ has the same number of atoms and coatoms; the following proposition appears in \cite{BHMPW-IH-module}.

\begin{proposition}
The following statements are equivalent:
\begin{itemize}
	\item ${}^p\!P_\cL(t) = 1$,
	\item ${}^p\! P_{\cL_F}(t) = 1$ for all $F \in \cL$,
	\item the lattice $\cL$ is modular, and for every rank $2$ interval $[F,G]\subset\cL$, we have $\#\{H \mid F < H < G\} \not\equiv 1 \!\!\mod p$.
\end{itemize}
\end{proposition}
For instance, if $\mathbb{K}$ is a finite field and $\cL$ is the poset of all vector subspaces of $\mathbb{K}^E$, then ${}^p\!P_\cL(t) = 1$ if and only if $p \neq \operatorname{char} \mathbb{K}$.


\bibliographystyle{siamplain}
\bibliography{symplectic}
\end{document}